\documentclass[12pt]{article}
\usepackage{enumerate,amsthm,amsfonts,amssymb,pstricks,hyperref, topcapt}%
\usepackage[leqno]{amsmath}
\usepackage{amscd}
\usepackage{amsfonts,rotating}
\usepackage{graphicx}%
\usepackage{fancyhdr}

\divide\oddsidemargin
by 2  

\makeatletter

\newcommand{\address}[2][]{%
  \ifx\@add@ress\@undefined\gdef\@add@ress{\par\par\bigskip}\AtEndDocument{\@add@ress}\fi
  \g@addto@macro\@add@ress{\bigskip\noindent{\small\scshape%
      \ifx#1\empty\else{\bfseries Address of #1:}\ \fi#2}\par\par}}

\renewenvironment{abstract}{\small\quotation\noindent
  {\bfseries \abstractname}}{\endquotation \par}
\newcommand{\footnotetextplain}[1]{\begingroup\def\@thefnmark{}%
  \long\def\@makefntext##1{\parindent 0pt\noindent ##1}\@footnotetext{#1}
  \endgroup}

\newcommand{\AMSsubjclass}[2]{\footnotetextplain{2000
   \emph{Mathematics Subject Classification:} Primary #1, Secondary #2.}}

\newcommand{\keywords}[1]{\footnotetextplain{\emph{Key words and phrases:} #1.}}

\makeatletter

\xdef\qedbuit{\qed}

\newcommand{\TeoremaAmbFinalMarcat}[1]{%
  \expandafter\gdef\csname end#1\endcsname{\qedbuit\@endtheorem}}

\theoremstyle{definition}
 \TeoremaAmbFinalMarcat{defi}
 \TeoremaAmbFinalMarcat{ex}
 \TeoremaAmbFinalMarcat{rem}
 \TeoremaAmbFinalMarcat{hrem}
 \TeoremaAmbFinalMarcat{conv}
 \TeoremaAmbFinalMarcat{prob}
 \TeoremaAmbFinalMarcat{conj}
 \TeoremaAmbFinalMarcat{conj}

 \newenvironment{proclama not emphasized}[1]{\par\vspace{\topsep}\noindent{\bf #1}}{\par\vspace{\topsep}}

\newenvironment{prooftext}[1] 
               {\noindent \textit{ #1}~ }     
               {\hfill\rule{2.5mm}{2.5mm} \vspace{\parskip} } 

\newcommand{\start}[2]{\begin{#1}\label{#2}}
\newcommand{\secc}[1]{Section~\ref{#1}}
\newcommand{\theoc}[1]{Theorem~\ref{#1}}

\newcommand{\figc}[1]{Figure~\ref{#1}}

\def\@enum@{\list{\csname label\@enumctr\endcsname}%
           {\usecounter{\@enumctr}\def\makelabel##1{\hss\llap{##1}}
           \itemsep=2pt\parsep=0pt\topsep=3pt plus 1pt minus 1 pt}}

\newenvironment{numlist}{\enumerate[(1)]}{\endenumerate}

\newcommand{\vs}{\vspace{.1in}}

\title{Impossible configurations for geodesics on negatively-curved surfaces}
\author{Anthony Phillips}
\date{\today}

\begin{document}
\maketitle

\begin{abstract}
  Hass and Scott's example of a 4-valent graph on the 3-punctured
  sphere that cannot be realized by geodesics in any metric of
  negative curvature is generalized to impossible configurations
  filling surfaces of genus $n$  with $k$ punctures for any $n$
  and $k$.

\end{abstract}

\AMSsubjclass{53C22}{57M50,30F99}

\section{Introduction}
By a configuration we mean a surface $S$ together with a 4-valent, connected
graph embedded in $S$. Going straight (neither right nor left) at
each intersection decomposes such a graph canonically into a
collection of closed curves (the {\em tracks})
intersecting themselves and each other
transversally.

\begin{figure}[htp]
\centering \includegraphics[width=1.5in]{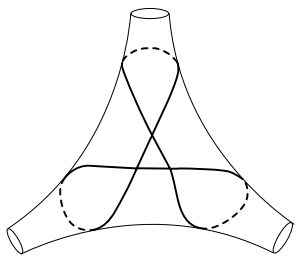} 
\caption{The Hass-Scott example.}
\label{hass-scott}
\end{figure}

A basic question is whether or not there is a negative-curvature metric on
$S$ such that the graph is isotopic to a collection of closed
geodesics intersecting transversally. We will say in this case that
the configuration can be {\em realized} by geodesics in a metric of
negative curvature.
It is an old but remarkable fact that the simple
configuration shown in \figc{hass-scott} cannot be so realized.
Joel Hass and Peter Scott
discovered this phenomenon in 1999 \cite{[HS]}. As              
they remark, their proof of non-realizability can be
by replaced by an argument, due to
Ian Agol, using the Gauss-Bonnet Theorem. 

In this paper Agol's argument is generalized to produce
  an infinite
  family of non-realizable configurations, the {\em polygonal impossible
    configurations}, including
  a surface
of genus $n$ with $k$ punctures for every $n\geq 2$ and $k\geq 0$.
In some sense polygonal impossible configurations
are all of the non-realizable examples that can be constructed
  using our general form of the argument.

\subsection{Preliminaries}\label{prelim}
The Hass-Scott example, where $S$ is the 3-punctured sphere,
is striking because at first inspection there
seems to be no reason for it not to be realizable. Locally, it
looks exactly like other, realizable configurations. In particular
\begin{itemize}
\item[C1.] The curve segments do not enclose embedded or immersed
  contractible 1-gons or 2-gons
  (``monogons'' or ``bigons''): these are in fact outlawed in
  geodesic configurations in negative curvature.
  \item[C2.] No curve represents a
 power ($>1$) in the
free homotopy group $\pi_1S$, nor do two distinct curves
  represent powers ($\geq 1$) of the same element of $\pi_1S$. Since
in negative curvature each free homotopy class contains a unique
geodesic, a power curve collapses to multiple tracings of
a single geodesic, and two homotopic curves collapse to the
same geodesic; in either case the initial configuration is
destroyed.
\end{itemize}

Furthermore, the Hass-Scott example {\it fills} its surface
in the sense that
\begin{itemize}
\item[C3.] The complement of the graph is a disjoint union of
  discs or singly punctured discs.
\end{itemize}
The configurations constructed here will satisfy the conditions
C1, C2 and C3.

\section{Polygonal impossible configurations}\label{pics}
\start{defi}{pic} A {\em polygonal impossible configuration}
$\mathcal{P}$
is an orientable, connected  2-dimensional 
surface
constructed as follows:
\begin{numlist}

\item Choose a number $N\geq 3$, which will be the number of vertices in
the configuration.
\item Choose a number $p$, with $N/4\leq p\leq N/3$, and
  $p$ polygons $A_1,\dots,
 A_p$ which together have $N$
corners
  (This is possible since $p\leq N/3$). At least one $A_i$ must
be a triangle: see the remark below.
\item Choose $q=N-2p$ even-sided polygons $B_1,\dots, B_q$ which
together have $2N$ corners. (This is possible since $p\geq N/4$ implies
$4q = 4N-8p \leq 2N$). Label the edges of the $B$-polygons alternately
active and inactive.
\item Identify an edge of one of the $A_i$ with every active edge of each
$B_j$, preserving orientations. Avoid forming a ring of squares: such
a ring would lead to two parallel tracks.
\end{numlist}
\end{defi}

\start{rem}{triangle}
At least one of the $A_i$ must be a triangle.
In fact, suppose first all the $A_i$ are squares;
then $N=4p, q=N-2p=2p$, and $2N/q=4$, so all the $B_i$ must also
be squares; and then the configuration $\mathcal{P}$
constructed by the algorithm will
be made up of one or more sets of parallel curves, contradicting C2 above.
On the other hand if all the $A_i$ have $\geq 4$
sides, and at least one has strictly more, then $N>4p$ contradicting
$p\geq N/4$.
\end{rem}

\start{theo}{impossible}
A polygonal impossible configuration 
cannot be given a metric of negative curvature
so that the set of curves defined
by its 1-skeleton is a set of geodesics.
\end{theo}
\begin{proof} Suppose such a metric exists.
Set $n_i$ to be the number of vertices of $A_i$, and $\alpha_{i,j}$,
$j=1,\dots, n_i$ to be the interior angle at the $j$th vertex of $A_i$.

Likewise set $m_i$ to be the number of vertices of $B_i$, and $\beta_{i,j}$,
$j=1,\dots, m_i$ to be the interior angle at the $j$th vertex of $B_i$.

Assume all the edges are geodesic arcs extending smoothly from
polygon to polygon, so that each of the $\alpha_{i,j}$ is
complementary to exactly two of the $\beta_{i,j}$.

The Gauss-Bonnet theorem \cite{[Wu]} gives

\begin{center}
$\begin{array}{ccc}
\alpha_{1,1} + \alpha_{1,2} + \cdots \alpha_{1,n_1}& <&(n_1-2)\pi\\
\dots& &\\

  \alpha_{p,1} +\alpha_{p,2} + \cdots \alpha_{p,n_p}& <& (n_p-2)\pi.
\end{array}                                                         $
\end{center}
Adding these equations,

$$\sum _{i=1}^p\sum_{j=1}^{n_i}\alpha_{i,j}< (N - 2p)\pi. ~~~~  (*)$$

Similarly, the sum of all the $\beta$s is strictly less than
$(2N - 2q)\pi$.
On the other hand each $\beta$ is $\pi - \alpha$ for some
$\alpha$, with each $\alpha$ occurring exactly twice.

So $$(2N-2q)\pi > \sum _{i=1}^q\sum_{j=1}^{m_i}\beta_{i,j}
= 2\sum _{i=1}^p\sum_{j=1}^{n_i}(\pi -\alpha_{i,j})
= 2N\pi -2\sum _{i=1}^p\sum_{j=1}^{n_i}\alpha_{i,j}$$
i.e. $\sum _{i=1}^p\sum_{j=1}^{n_i}\alpha_{i,j} > q\pi$.
Since by the construction $q=N-2p$, this inequality contradicts $(*)$.
\end{proof}

\section{The genus of a polygonal impossible configuration; minimal and unicursal configurations}\label{genus}
A polygonal configuration $\mathcal{P}$ created from $N$,
$A_1, \dots, A_p$, $B_1, \dots, B_q$ as above comes with a surface in
which it is naturally embedded:  
the
$N$ inactive edges of the $B$s are grouped by the identifications in
step 4 into a collection $\gamma_1, \dots, \gamma_r$ of
closed curves;
adding a disc $D_i$ along each $\gamma_i$ creates a closed
orientable surface $S_{\mathcal P}$.

The surface $S_{\mathcal{P}}$ has  Euler characteristic
$\chi=N-2N+(p+q+r) = -p+r~$, and genus
$$g_{\mathcal{P}} =
\frac{1}{2}(2-\chi) = \frac{1}{2}(2-r+p)~~~~~(**).$$
We can take $g_{\mathcal{P}}$
as the genus of $\mathcal{P}$; this matches the usual definition of
the genus of a graph as the genus of the simplest surface on which it
can be embedded so that its complement is topologically a set of discs.\vs

The surface $S_{\mathcal{P}}$ does not necessarily admit
a metric of negative curvature.

{\it Preliminary punctures.} 
To start, some of the discs $D_1, \dots, D_r$ may be bounded by a single
or exactly two curve segments. To satisfy condition C1 these discs
must be punctured.

{\it Further necessary punctures.} The next steps depend on $g_{\mathcal{P}}$.
\begin{itemize}
\item $g_{\mathcal{P}} =0$. Since $p\geq 1$ the equation $(**)$ implies
  $r\geq 3$. Puncturing three of $D_1, \dots, D_r$ if still necessary will give the
  3-punctured sphere, a surface   
  admitting a metric of negative curvature. Note that a sphere with one
  or two punctures admits such a metric, but in the first case there are
  no geodesics (so every configuration is impossible), and in the second case
  the only possible configuration is a circle.
  
\item $g_{\mathcal{P}} =1$. Here $(**)$	implies $r\geq 1$;
  puncturing one of $D_1, \dots, D_r$ if still necessary gives the punctured torus,
  a surface	
  admitting a metric of	negative curvature.

\item $g_{\mathcal{P}} \geq 2$. In this case the surface admits
  a metric of negative curvature.

\end{itemize}

\start{defi}{minimal}
A polygonal configuration on a surface a surface of genus $g$
with $k$ punctures is {\em minimal} if
it fills the surface, and has the smallest possible number of vertices
for such a configuration.
\end{defi}

This minimal number depends on $g$ and $k$.
\begin{itemize}
  \item If $k=0$, since $r\geq 1$
equation $(**)$ implies that $p$, the number of $A$-polygons,
  must satisfy $p\geq 2g-1$. Since each $A$-polygon is at least
  a triangle, the number $N$ of vertices for a polygonal impossible
  configuration 
  filling the unpunctured
  surface of genus $g$ must satisfy $N\geq 6g-3$.
\item Otherwise, since $r\geq k$ we have $p\geq 2g+k-2$
  and  $N\geq 6g+3k-6$.
\end{itemize}

  \start{defi}{unicursal}
A configuration is {\em unicursal} if it has
exactly one track, i.e., as described in the introduction,
if it can be traversed by a single curve.
\end{defi}

\section{Examples}\label{examples}

\begin{figure}[htp]
\centering \includegraphics[width=3in]{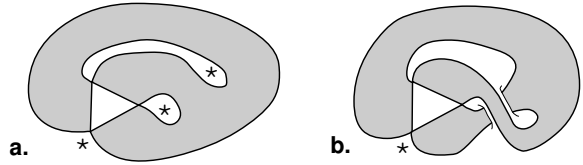} 
\caption{$N=3$. (Asterisks represent punctures).
{\bf a.} Here the identification gives
$r=3$ and $g_{\mathcal{P}}=0$, yielding  a unicursal
impossible configuration on the 3-punctured sphere (the
configuration exhibited by Hass and Scott).
{\bf b.} Otherwise the identification gives $r=1$ and $g_{\mathcal{P}}=1$;  hence
a graph on the torus which becomes an impossible (3-track)
configuration
on the punctured torus.} 
\label{N=3}
\end{figure}

\begin{figure}[htp]
\centering \includegraphics[width=5in]{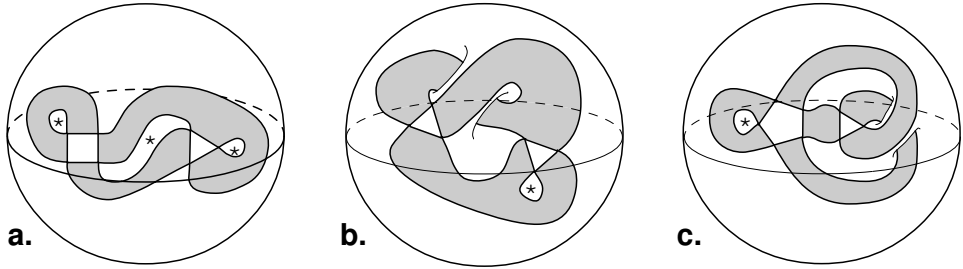} 
\caption{{\bf a.} a unicursal
  configuration on the 3-punctured sphere,  {\bf b.} and {\bf c.},
  2-track configurations on the
  punctured torus. }
  
\label{pentagon}
\end{figure}

The smallest possible $N$ is $N=3$. Here $p$ and $q$  must equal 1,
with $A_1$ a triangle and $B_1$ a hexagon.
There are two ways to make the identification, with different results,
as shown in \figc{N=3}. Additional simple examples are shown in \figc{pentagon}.

\section{Unicursal, filling configurations on surfaces of higher genus}

\subsection{Surfaces of                
  genus $\geq 2$}
\start{prop}{fill-smooth}
There exists a minimal unicursal impossible polygonal
configuration filling the surface of genus $n$.
\end{prop}

\begin{proof} We begin with genus $2$.
  \vs
  \newpage
  
  \begin{figure}[htb]
\centering \includegraphics[width =1.5in]{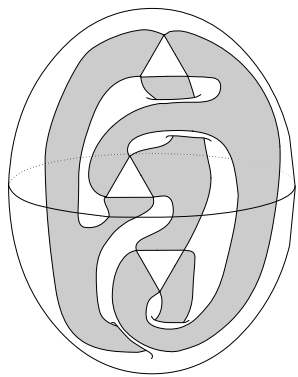}
\caption{The 1-track impossible configuration $\mathcal{P}_2$ on the genus-2 surface with
  $N=9, p=3, q=3$ (3 hexagons), $r=1$.}
\label{genus-2}
\end{figure}

This construction can be extended to give an impossible configuration $\mathcal{P}_n$
filling
the surface of genus $n$, for each $n\geq 2$.

\begin{figure}[htb]
  \centering \includegraphics[width =3in]{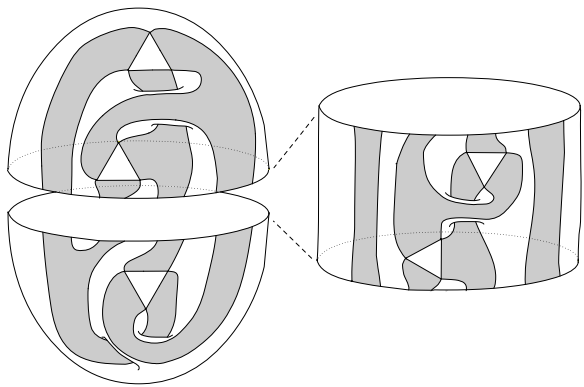}

    \caption{The configuration $\mathcal{P}_n$ is obatined by
  interpolating into the $\mathcal{P}_2$ configuration $n-2$ copies of the partial configuration shown
  on the right. It 
  has $p=2n-1$ triangles 
matched with $q=2n-1$ hexagons. As can be checked, it has $r=1$  and
therefore Euler characteristic $-2n+2$ and genus $n$. }
\label{genus-n}
\end{figure}

Note that the  configuration $\mathcal{P}_n$
is unicursal for all $n\geq 2$. Also, since $\mathcal{P}_n$
has $2n-1$ triangles and hence $6n-3$ vertices, it is minimal.
\end{proof}

\subsection{Surfaces with punctures}

The Hass-Scott example is minimal and unicursal. Generalizing it
to surfaces with more punctures, using polygonal impossible
configurations, can be done preserving both of
these properties if the number of punctures is odd,
but only one or the other if it is even.

\start{prop} {punctures}
The configuration $\mathcal{P}_n$ can be extended to become
an impossible configuration filling the surface of genus $n$
with $k$ punctures, for any $k\geq 0$. If $k$ is odd, the
new configuration can be minimal and unicursal. If $k$ is
even, it can be minimal or unicursal but not both.
\end{prop}

\begin{proof}
Initially $\mathcal{P}_n$ has one complementary region, a disc. 
This disc can be punctured, yielding $k=1$. Splicing in $m$ copies
of partial configuration {\bf c} (Fig. \ref{punctures}) gives a
minimal, unicursal polygonal configuration with the same genus
and $2m+1$ punctures. On the other hand as remarked in \secc{genus}
a minimal configuration of genus $n$ with $k$ punctures has
$N=6n+3k-6$ vertices; if $k$ is even, so is $N$; and by
\theoc{parity} in the Appendix,
such a configuration cannot be unicursal.

A unicursal configuration with
genus $n$, with $2m$ punctures and one extra vertex can be obtained from
$\mathcal{P}_n$ by splicing in
$m-1$ copies of {\bf c} and one of {\bf b}. A minimal configuration
with genus $n$, with $2m$ punctures and two tracks can be obtained by splicing in
$m-1$ copies of	{\bf c} and one of {\bf a}.
  \begin{figure}[htb]
\centering \includegraphics[width =3in]{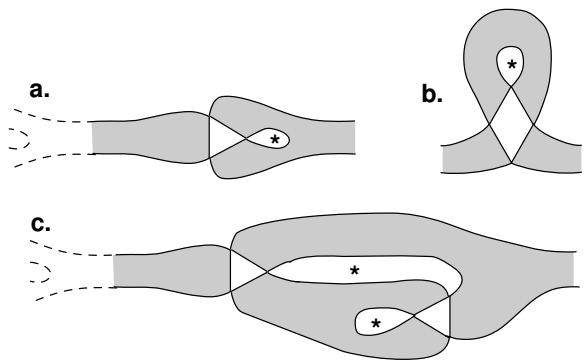}
\caption{The $*$ represents a puncture. These partial configurations
  can be spliced into $\mathcal{P}_n$ to increase the
  number of punctures by one ({\bf a}, {\bf b}) or two ({\bf c}).
  Note that partial configurations {\bf a}
  and {\bf c} 
  must be spliced onto an arm coming {\it away} from a hexagon, as
  shown, in order
  to avoid producing a 2-gon.
}
\label{punctures}
\end{figure}
\end{proof}

\section{Appendix: Unicursal configurations require an odd number of vertices}
The examples in Figs. \ref{N=3}, \ref{pentagon} and \ref{genus-2}
suggest the following statement.

\start{theo} {parity}
The number of tracks of a
polygonal impossible configuration is congruent mod 2 to the
number of vertices. In particular, such a configuration can only be
unicursal if the number of vertices is odd.
\end{theo}

{\bf Preliminaries for the proof.}

\begin{numlist}
  \item Taking the planar polygons
$A_i$ and $B_i$ as in the construction of the
configuration, give each one the standard (counterclockwise)
orientation. Then give the segments of the configuration
their inherited orientation, except segments shared by an
$A$ and a $B$ keep their
$A$-orientation. With this convention, the track-segments
  at each intersection are coherently oriented, and give  a well-defined
orientation on each track (\figc{segment-orientations}).

\begin{figure}[htb]
\centering \includegraphics[width=3in]{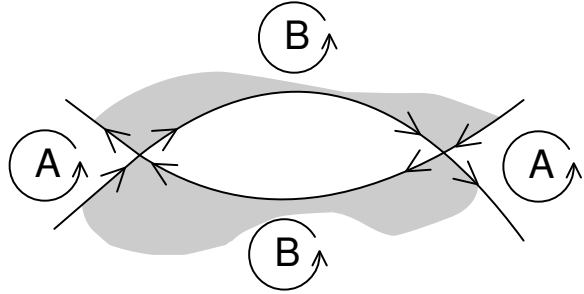}
\caption{The $A$ and $B$ orientations are adjusted to give
a well-defined
orientation on each track.}
\label{segment-orientations}
\end{figure}

\item Project the configuration into the plane, and consider
it as a collection of oriented immersed curves.
The configuration depends only on the
nature of the polygons $A_i$ and $B_i$ and the way they are connected.
In particular, its projection can be displayed so that the
$A_i$ and $B_i$ appear as in \figc{projection}.
\begin{figure}[htb]
\centering \includegraphics[width=5in]{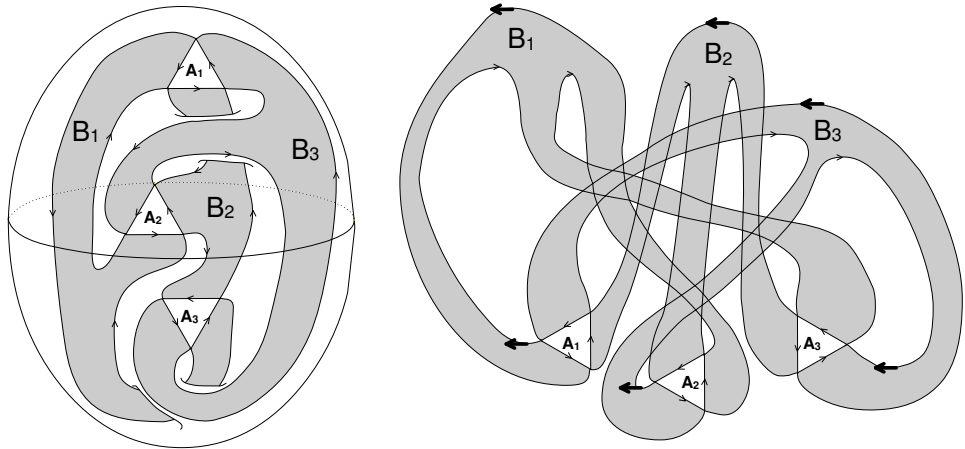}
\caption{The projection of a
  polygonal impossible configuration can be displayed
with all the $B$ polygons at the top, all the $A$ polygons at
the bottom, and so that the only horizontal tangents appear at
the top and at the bottom. This figure shows the configuration
from \figc{genus-2}, oriented as above, with its display in this form.
Bold arrows: the unit tangent vector is equal to $(-1,0)$. }
\label{projection}
\end{figure}
\end{numlist}

\begin{prooftext}{Proof of \theoc{parity}}
  \begin{numlist}
\item We show that the sum of the rotation numbers of the projected
complex of curves
is even. For each curve, the rotation number can be defined as the degree
of the Gauss map, which takes a
parameter value to the corresponding unit tangent vector,
considered as a point on the unit circle. The degree of a
smooth map is equal modulo 2 to the number of inverse images of
a regular value \cite{[M]}. For a regular value we choose $(-1,0)$, the horizontal
unit vector pointing left. \vs

First, inspection of \figc{projection} shows that each
$B$ polygon contributes exactly one to the count of inverse images
of $(-1,0)$. With notation from the definition of polygonal impossible
configuration, the contribution of the $B$-polygons is $q$.
\vs

Next we will show that $(***)$ the $n_i$-gon  $A_i$
contributes $n_i-2$ to this count. It will follow that the total
contribution of the $A$-polygons is $\sum_{i=1}^p (n_i - 2) = N-2p$.
Since $N-2p=q$, adding in the contribution of the $B$s gives as total the even number $2q$.
\begin{figure}[htb]
\centering \includegraphics[width=3in]{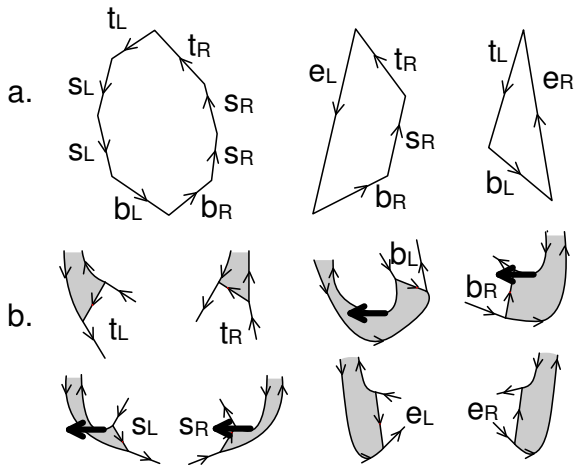}
\caption{{\bf a.} The eight possible relative positions of an edge
  in a general-position oriented polygon:
  ${\bf t}_L, {\bf t}_R$ top left and right,
${\bf s}_L, {\bf s}_R$ side left and right, ${\bf b}_L, {\bf b}_R$
botton left and right and ${\bf e}_L, {\bf e}_R$ top-to-bottom left
and right. 
{\bf b.} The edges  ${\bf s}_L, {\bf s}_R$, ${\bf b}_L, {\bf b}_R$ each contribute 1
to the count of inverse images of $(-1,0)$ under the Gauss map; the
other edges do not contribute.}
\label{corners}
\end{figure}
\vs

 To prove $(***)$ note (\figc{corners}) that an oriented convex $n$-gon in
general position (no horizontal sides) must
have  one of these three configurations:
 \begin{numlist}
 \item ${\bf t}_L + {\bf t}_R + {\bf b}_L + {\bf b}_R +
  (n-4)\{ {\bf s}_L, {\bf s}_R\}$
  \item ${\bf e}_L + {\bf t}_R + {\bf b}_R + (n-3){\bf s}_R$
  \item ${\bf e}_R + {\bf t}_L + {\bf b}_L + (n-3){\bf s}_L$.
  \end{numlist}
  In each of these cases, the number of inverse images of $-x$ is $n-2$.

\item The number of self-intersection points of an immersed oriented curve
in the plane, counted mod 2, is one less than its rotation number.
(Because the self-intersection number mod 2 is a regular homotopy
invariant \cite{[W2]}, and because any
curve with rotation number $n$ is regularly homotopic to $n$ turns of a
spiral, with the endpoints joined \cite{[W1]}:
a curve with $|n|-1$ intersection points).
So a curve with even rotation number must have an odd number
of self-intersection points.
\item Let $\gamma_1, \dots, \gamma_k$ be the $k$ tracks of the path
through our configuration. We know that the sum of their winding numbers
is even, so an {\em even} number $\ell$ of them have odd winding
number; these $\ell$ tracks each have even self-intersection number.
The other $k-\ell$ tracks have even winding number and therefore
odd self-intersection number. The sum of their self-intersection numbers
is therefore congruent to $k-\ell$ and therefore to $k$, since
$\ell$ is even; it follows that the sum of all the self-intersection numbers
of the $\gamma_i$ is congruent mod 2 to $k$.
\item Finally, the self-intersection points of the path through
the configuration, drawn as in \figc{projection},  are of three types: those
coming from the self-intersection numbers of $\gamma_i$ for $i=1,\dots,k$,
those coming from intersections between $\gamma_i$ and $\gamma_j$ for
$i\neq j$, and those
coming from the intersections of the descending arms of the $B$-polygons.
The second and third types come in pairs. So the total number of
self-intersections is congruent mod 2 to $k$, the number of tracks; and this must also hold
for the number of those of the first two types, which is the number
  of vertices of the configuration.
\end{numlist}    
  
\end{prooftext}

\bibliographystyle{amsalpha}

{\sc

\noindent
Department of Mathematics,
Stony Brook University,
Stony Brook NY 11794
}

\noindent
\emph{E-mail}{:\;\;}\texttt{tony@math.stonybrook.edu}

\end{document}